\documentclass[10pt]{article}
 
\usepackage[latin1]{inputenc}
\usepackage{relsize}
\usepackage{tikz} 
\usepackage{framed} 
\usepackage{url, amsmath,enumerate,fancyhdr,amssymb, amsthm, 
pstricks, epsf, lscape, ifpdf, graphicx, floatrow, multirow} 
\usepackage{algorithm}
\usepackage{algorithmicx} 
\usepackage{algpseudocode}
\usepackage{enumitem}

\usepackage{url, amsmath,enumerate,fancyhdr,amssymb, amsthm, 
	pstricks, epsf, lscape, ifpdf, graphicx, float}
\usepackage[pagebackref=true, colorlinks=true, linkcolor=blue, citecolor=blue, linkcolor=blue, anchorcolor=red, urlcolor=blue]{hyperref}
\usepackage{algorithm}
\usepackage[margin=3cm]{geometry}
\usepackage{mathtools}
\usepackage{multicol}
\usepackage{tabularx}
\newcolumntype{C}{>{\centering\arraybackslash}X}
\usepackage{amsmath}
\makeatletter
\newcommand{\leqnomode}{\tagsleft@true}
\newcommand{\reqnomode}{\tagsleft@false}
\makeatother

\newcommand{\inProd}[2]{\langle #1 , #2 \rangle }

\usetikzlibrary{shapes.misc,shadows}

\algnewcommand\algorithmicinput{\textbf{Input:}}
\algnewcommand\INPUT{\item[\algorithmicinput]}
\algnewcommand\algorithmicinitialization{\textbf{Initialization:}}
\algnewcommand\INITIALIZATION{\item[\algorithmicinitialization]}

\usepackage{longtable,tabu}

\hypersetup{
	colorlinks,
	linkcolor={red!50!black},
	citecolor={blue!50!black},
	urlcolor={blue!80!black}
}

\iffalse 
\newcommand{\myblue}[1]{\textcolor{blue}{#1}} 
\else
\newcommand{\myblue}[1]{#1} 
\fi

\usepackage{tabularx}

\renewcommand{\theequation}{\thesection.\arabic{equation}}
\newcommand{\pha}{\phantom}

\mathchardef\mhyphen="2D

\newcommand{\val}{\operatorname{val}}

\newtheorem{Definition}{Definition}

\newtheorem{Example}{Example}
\newtheorem{Proposition}{Proposition}
\newtheorem{Lemma}{Lemma}
\newtheorem{Theorem}{Theorem}
\newtheorem{Corollary}{Corollary}
\newtheorem{Remark}{Remark}
\newtheorem{Assumption}{Assumption}

\newcommand{\A}{\mathcal A}

\newcommand{\ri}{\operatorname{ri}}

\newcommand{\trace}{\operatorname{tr}}

\newcommand{\eps}{\epsilon} 
\newcommand{\bpx}{\begin{pmatrix}}
\newcommand{\epx}{\end{pmatrix}}
\newcommand{\bbx}{\begin{bmatrix}}
\newcommand{\ebx}{\end{bmatrix}}

\newcommand{\bdef}{\begin{Definition}} 
\newcommand{\commentout}[1]{}
\newcommand{\co}[1]{}

\newcommand{\nin}{\noindent}
\newcommand{\ti}{\times}
 
\newcommand{\pf}[1]{\vspace{.35cm} \nin {\bf Proof {#1} }}

\newcommand{\norm}[1]{\lVert  #1  \rVert} 
\newcommand{\sym}[1]{{\cal S}^{#1}}
\newcommand{\psd}[1]{{\cal S}_+^{#1}}
\newcommand{\mymatrix}[1]{\rad{n \times n}}

\newcommand{\rad}[1]{\mathbb{R}^{#1}}

\newcommand{\symn}{\sym{n}}
\newcommand{\psdn}{\psd{n}}

\newcommand{\R}{ {\mathcal R} }

\newcommand{\T}{{\tan}}

\newcommand{\la}{\langle}
\newcommand{\ra}{\rangle}

\newcommand{\dist}{\operatorname{dist}}

\newcommand{\beq}{\begin{equation}}
\newcommand{\eeq}{\end{equation}}
\newcommand{\beqa}{\begin{eqnarray}}
\newcommand{\eeqa}{\end{eqnarray}}
\newcommand{\ba}{\begin{array}}
\newcommand{\ena}{\end{array}}
\newcommand{\bac}{\begin{array}{ccccccccccc}}
\newcommand{\eac}{\end{array}}
\newcommand{\bprop}{\begin{Proposition}}
\newcommand{\eprop}{\end{Proposition}}

\newcommand{\beqast}{\begin{eqnarray*}}
\newcommand{\eeqast}{\end{eqnarray*}}
\newcommand{\benum}{\begin{enumerate}}
\newcommand{\eenum}{\end{enumerate}}
\newcommand{\bit}{\begin{itemize}}
\newcommand{\eit}{\end{itemize}}
\newcommand{\bth}{\begin{Theorem}}
\newcommand{\enth}{\end{Theorem}}
\newcommand{\ble}{\begin{Lemma}}
\newcommand{\ele}{\end{Lemma}}
\newcommand{\bex}{\begin{Example}}
\newcommand{\eex}{\end{Example}}
\newcommand{\bcor}{\begin{Corollary}}
\newcommand{\ecor}{\end{Corollary}}
\newcommand{\brem}{\begin{Remark}}
\newcommand{\erem}{\end{Remark}}
\newcommand{\bass}{\begin{Assumption}}
\newcommand{\eass}{\end{Assumption}}

\renewcommand{\theequation}{\thesection.\arabic{equation}}

\newcommand{\bsmx}{\begin{small} \begin{pmatrix}}
\newcommand{\esmx}{\end{pmatrix} \end{small}}

\newcommand\bovermat[2]{%
	\makebox[0pt][l]{$\smash{\overbrace{\phantom{%
					\begin{matrix}#2\end{matrix}}}^{\text{#1}}}$}#2}

\usepackage{mathtools}
\makeatletter
\renewcommand*\env@matrix[1][*\c@MaxMatrixCols c]{%
	\hskip -\arraycolsep
	\let\@ifnextchar\new@ifnextchar
	\array{#1}}
\makeatother

\usepackage{tikz,array,calc, xcolor}
\usetikzlibrary{decorations.pathreplacing}

\begin{document}


\title{A Simplified Treatment of Ramana's Exact Dual for Semidefinite Programming}  
\author{Bruno  F. Louren\c{c}o\thanks{Department of Statistical Inference and Mathematics, Institute of Statistical Mathematics, Japan. (\texttt{bruno@ism.ac.jp})} 
	\and 
	G\'{a}bor Pataki\thanks{Department of Statistics and Operations Research, University of North Carolina at Chapel Hill, USA. (\texttt{gabor@unc.edu})}
}

\maketitle 

\begin{abstract}
{\color{black} 
		In semidefinite programming  the dual may fail to attain 
		its optimal value and there could  be a duality gap, i.e., the primal and dual 
		optimal values may  differ. In a striking paper, Ramana~ \cite{Ramana:97} proposed a polynomial size extended dual that  does not have these deficiencies and yields a number of fundamental  results  in complexity theory. In this work we walk the reader through a concise and self-contained derivation of Ramana's dual,
		relying mostly on elementary linear algebra.
		}
\end{abstract}



{\em OR/MS subject classification:} Primary: convexity; secondary: programming-nonlinear-theory

\section{Introduction} \label{sect-intro}

Consider the primal-dual pair of semidefinite programs (SDPs) 
	\begin{center}
	\begin{minipage}{0.5\linewidth}
		\leqnomode
		\begin{equation}\label{p}
			\begin{split}
	        \sup & \,\, \sum_{i=1}^m  c_i x_i \\
		s.t.   & \,\,  \sum_{i=1}^m  x_i A_i \preceq B \\
		\end{split}\tag{P}
		\end{equation}
	\end{minipage}%
		\begin{minipage}{0.5\linewidth}
		\begin{equation}\label{d}
		\begin{split}
		\inf  & \,\, \la B,  Y  \ra  \\
		s.t. & \,\, \la A_i, Y \ra \, = \, c_i \, (i=1, \dots, m) \\
		& \,\, Y \succeq 0
		\end{split}\tag{D}
		\end{equation}
	\end{minipage}
\end{center}
where  $A_1, \dots, A_m, \,$ and $B$ are $n \times n$ 
symmetric matrices and  $c_1, \dots, c_m$ are scalars. 
For  symmetric matrices $S$ and $T$  we write 
$S \preceq T$ to say that $T - S$ is positive semidefinite (psd) 
and we write $\la T, S \ra :=  {\rm trace}(TS)$ to denote their inner product.

 SDPs are an elegant generalization of linear programming, and they appear in a broad range of application areas. 
 However,  the duality theory of SDPs is much less satisfactory than 
 that of  linear programming. \myblue{On the one hand, the optimal value of \eqref{p} is bounded from above by the optimal value of \eqref{d}.
 On the other, } \eqref{p} and \eqref{d} may not have optimal solutions, i.e., 
 the supremum may not be a maximum, and the infimum may not be a minimum.
 Even worse,  the optimal values of \eqref{p} and \eqref{d} may differ.

\bex \label{ex-main}
In the following classical  pathological SDP 
	\beq \label{ex1-problem}
	\begin{array}{rl}
		\sup &   2 x_1 \\
		s.t. & x_1 \bpx 0 & 1 \\ 1 & 0 \epx \preceq \bpx 1 & 0 \\ 0 & 0 \epx
	\end{array} 
	\eeq
	the constraint is equivalent to 
	        $
	        \bigl( \begin{smallmatrix} 1         & - x_1 \\
	        - x_1    & 0 \end{smallmatrix} \bigr) \succeq 0,
	        $
	         so the only feasible solution is $x_1 = 0.$ 
	         
	         The dual, with a variable matrix 
	        $Y = (y_{ij}), \,$
	        is 
	       	\beq \label{ex1-problem-dual} 
		\ba{rrcl}
		\inf & y_{11} \\
		s.t. & y_{12}  & = &  1 \\
		      & Y & \succeq & 0,
		\ena
		\eeq
		wherein any $Y_\epsilon :=  \bigl(  \begin{smallmatrix} \eps & 1 \\ 1 & 1/\eps \end{smallmatrix} \bigr)$ 
		with $\epsilon > 0$ is feasible. So we conclude that the infimum of \eqref{ex1-problem-dual} 
		is $0.$ 
			
				However, any $Y$ with $y_{11} = 0$ is not feasible  in \eqref{ex1-problem-dual}, so its infimum  is not attained. 
\eex 
		In a striking paper, Ramana \cite{Ramana:97} constructed a new dual problem that fixes most of the issues of the classical SDP dual. Ramana's dual has the following  attractive traits:
		\begin{enumerate}
			\item it does not assume anything about 
			\eqref{p}, other than it is feasible; 
			\item it attains its optimal value, when that value  is finite; 
			\item  its optimal value is the same as that of \eqref{p}, so there is no duality gap; 
				\item it yields important complexity implications. Among other things,  
			it proves that deciding feasibility of SDPs in the Turing model is not NP-complete, unless 
			NP = co-NP, which is an unlikely scenario according to most experts.
				\end{enumerate} 
					Ramana's dual  sparked great excitement in the SDP community and inspired many followup papers. 
		Ramana, Tun\c{c}el and Wolkowicz \cite{RaTuWo:97}  connected it to  the facial reduction algorithm of Borwein and Wolkowicz  \cite{BorWolk:81B, BorWolk:81};
		 Luo, Sturm, and Zhang \cite{LuoSturmZhang:97} gave a different proof of its correctness; 
	 Ramana and Freund \cite{RaFreund:96} showed that it has zero duality gap with its 
		usual dual; and Klep and Schweighofer \cite{KlepSchw:12} 
		constructed a dual with similar properties, which relies on 
		machinery from real algebraic geometry.   
		Ramana's work 
		is often cited in surveys and books:  see for example Drusvyatskiy and Wolkowicz \cite{drusvyatskiy2017many}, 
		DeKlerk \cite{de2006aspects}, Vandenberghe and Boyd  \cite{vandenberghe1996semidefinite}, Nemirovski \cite{nemirovski2007advances}, and 
	Laurent and Rendl 	\cite{laurent2005semidefinite}. 
	It was used by DeKlerk, Roos and Terlaky \cite{de2000self}  in 	self-dual embeddings. 
	 It is often mentioned in the discrete mathematics and theoretical computer science literature,  see for example, Lov{\'a}sz~\cite{lovasz2019graphs} and O' Donnell~\cite{o2017sos}.  
	 
			Ramana's dual has been generalized in a number of directions: to conic linear programs over so called 
		{\em nice cones} \cite[Corollary 1]{Pataki:13}, and   even to arbitrary conic linear programs, which  have a Ramana type dual 
	  \cite[Theorem 	2]{liu2018exact}.  It has greatly inspired the authors to examine why the pathologies arise 
	  in the first place: see for example \cite{lourencco2016structural} and \cite{pataki2019characterizing}. 
	  
	The known derivations of Ramana's dual rely on convex analysis, namely on the technique of facial reduction. Facial reduction originated in 
	 in the eighties  \cite{BorWolk:81B, BorWolk:81}, then simplified variants were proposed by Waki and Muramatsu \cite{WakiMura:12} and the second author of this note
	\cite{Pataki:13}. For a recent survey of facial reduction and its applications, we refer to  \cite{drusvyatskiy2017many}.  On the other hand,
	the related dual of Klep and Schweighofer \cite{KlepSchw:12}  employs algebraic geometry. These are two complementary approaches, tailored to readers trained 
	either in convex analysis, or algebraic geometry. 
		
		 In this work we give a short and elementary derivation of Ramana's dual that  we hope  will 
	appeal to all audiences. 
	
	To set the stage, we define the operator $\A$ and its adjoint $\A^*$ as 
\begin{equation} \nonumber 
\A x := \sum_{i=1}^m x_i A_i, \, \A^* Y = ( \inProd{A_1}{Y}, \dots,  \inProd{A_m}{Y} )^\top, 
\end{equation}
where $x \in \rad{m}$ and $Y$ is an $n \times n$ symmetric matrix.

In what follows, we assume that the primal \eqref{p} is feasible, and we denote by 
$\val()$ the optimal value of an optimization problem. 
We use the common convention that the optimal value of an ``inf" problem is $+ \infty$ exactly when
it is infeasible. 
We denote by $\sym{n}$ the set of $n \times n$ symmetric matrices, 
and by $\psdn$  the set of symmetric psd matrices. For a matrix 
$M$ (symmetric or not) $\R(M)$ stands for its range space. 
\begin{Theorem}\label{theo:main}
	Consider the optimization problem 
		\begin{equation}
		\left.
		\begin{minipage}{0.8\displaywidth}
		\addtocounter{equation}{-1}
		\renewcommand{\theHequation}{x\theequation}
		\vspace*{-\baselineskip}
		\begin{align}
		{\inf} &  \,\,\, \langle B, U_{n+1} + V_{n+1} \rangle  &     &             &                             & \nonumber \\ 
		s.t.    & \,\,\, \A^*(U_{n+1} + V_{n+1})                          & = & \,\, c    &                             & \nonumber \\ 
		& \A^* (U_{i} + V_{i})                                           & = & \,\, 0    & \!\!\!\!\!\!\!\!\! i = 1,\ldots, n     & \label{eqn-ai0}  \\ 
		& \langle B, U_{i} + V_{i} \rangle                     & = & \,\, 0    & \!\!\!\!\!\!\!\!\!  i = 1,\ldots, n     & \label{eqn-b0} \\
		& V_i \in \T(U_{i-1})    &     &             & \!\!\! i = 1,\ldots, n+1 &  \label{eqn-tan} \\
		& U_{i} \in \psd{n}                                                &     &             &  \!\!\! i = 1,\ldots, n+1  &   \label{eqn-Uipsd}   \\
		& U_{0} = V_0 = 0.                                               &    &             &                               & \nonumber 
		\end{align} 
		\end{minipage}
		\,\,\, \right\}
		\label{d-ramana}  \tag{D$_{\text{Ram}}$} 
		\end{equation}
			\addtocounter{equation}{1}
		
			Here for $U \in \psd{n}$ the set $\T(U)$ is defined as 
			\begin{equation} \label{eqn-tan:2}
			\T(U) \, = \, \{  \,  W + W^\top \, | \, W \in \rad{n \times n}, \, \R(W) \subseteq \R(U) \,   \},
			\end{equation}
		and called the {\em tangent space of} $\psd{n}$ at $U.$
		
		Then 	
		$$\val \eqref{p} = \val \eqref{d-ramana}, \,$$ and $\val \eqref{d-ramana}$ is attained when finite. 
		\qed
\end{Theorem}
We call \eqref{d-ramana} the {\em Ramana dual} of \eqref{p}. 

\co{ Alternative way to write Ramana dual:
	\begin{equation}
	\left.
	\begin{minipage}{0.8\displaywidth}
		\addtocounter{equation}{-1}
		\renewcommand{\theHequation}{x\theequation}
		\vspace*{-\baselineskip}
		\begin{equation} \label{d-ramana-test}  \tag{D$_{\text{Ram}}$} 
		\begin{array}{rrclll}
			{\inf} &  \,\,\, \langle B, Y  \rangle  &     &             &                             & \nonumber \\ 
			s.t.    & \,\,\, \A^* Y                               & = & \,\, c    &                             & \nonumber \\ 
			         & Y                                              & \in & \psd{n} + \tan(U_n) \\
			& \A^* (U_{i} + V_{i})                                           & = & \,\, 0    & \!\!\!\!\!\!\!\!\! i = 1,\ldots, n     &   \\ 
			& \langle B, U_{i} + V_{i} \rangle                     & = & \,\, 0    & \!\!\!\!\!\!\!\!\!  i = 1,\ldots, n     &  \\
			& U_{i} \in \psd{n},  V_i                          & \in     &   \T(U_{i-1})           &  \!\!\!\!\!\!\!\! i = 1,\ldots, n  &     \\
			& U_{0} = V_0 & = & 0,                                               &                                   & \nonumber 
		\end{array} 
	\end{equation} 	
	\end{minipage}
	\,\,\, \right\}
\end{equation}
}
Ramana's dual at first may look mysterious. The reader may also object
 that, because of  the tangent space 
 constraint in \eqref{eqn-tan},  it is not even an SDP! We will take care of the latter issue soon, but we
 first explain the makeup of \eqref{d-ramana}. 
 Variables $U_0$ and $V_0$ are  included in it only for convenience, and variable $V_1$ is  always zero since $U_0 = 0.$

It is straightforward that 
	\begin{equation} \label{eqn-dram-d} 
	\val \eqref{d-ramana} \leq \val \eqref{d}
\end{equation}
	holds. Indeed, we can construct a feasible solution of \eqref{d-ramana} from any feasible solution $Y$ of \eqref{d}: 
	we can set $U_{n+1} := Y$ and the other $U_i$ and all $V_i$ to zero. 
	 This construction works,  since 
	the zero matrix is in the tangent space of $\psd{n}$ at any psd matrix: see the definition in \eqref{eqn-tan:2}.
	
	 However, in general 
	$U_1, \dots, U_n$ and 
	the $V_i$ will not be all zero. Loosely speaking, they ensure that $U_{n+1} + V_{n+1}$ lives in the ``right" superset of $\psdn$ to ensure zero duality gap and attainment in \eqref{d-ramana}. 
We illustrate this fact 
	next. 

\bex  \label{ex1-continued}    (Example \ref{ex-main} continued) The Ramana dual of 
\eqref{ex1-problem} has an attained $0$ optimal value, as evidenced by the  solution 
\begin{equation} \label{eqn-UiVi-ex1} 
\begin{split}
U_1 \, = \, V_1 = 0, \, U_2 \, = \, \bpx 0 & 0 \\ 0 & 1 \epx, \,  V_2 \, = \, 0, \,\, U_3  \, =  \, 0, \,\,  V_3  \, =  \,  \bpx 0 & 1 \\ 1 & 0 \epx. 
\end{split} 
\end{equation}
\myblue{Here $V_3 \in \tan(U_2)$ because $V_3$  can be written as 
$$
V_3 \, = \, \underbrace{\bpx 0 & 0 \\ 1 & 0 \epx}_{W}    +   \underbrace{\bpx 0 & 1 \\ 0 & 0  \epx}_{W^\top}, 
$$
and $\R(W) \subseteq \R(U_2).$ } Note that $U_3 + V_3$ is not psd.  
\eex
In this example, as well as in later examples, we only list the $U_i$ and $V_i$ for $i \geq 1$ (since $U_0 = V_0 = 0$).

 The next theorem, which uses ideas from Lemma 2.1 in \cite{RaTuWo:97}, 
   shows that \eqref{d-ramana} can be turned into a bona fide SDP:
\begin{Theorem} \label{theorem-repr} 
	Suppose $U \in \psd{n}.$ Then 
	\begin{equation} \label{eqn-UVW} 
	\T(U) \, = \, \biggl\{  W + W^\top \, \bigl| \, W \in \rad{n \times n}, \, \begin{pmatrix}
	U   & W \\
	W^\top & \beta I
	\end{pmatrix}  \succeq 0 \, \text{for some} \, \beta \in \rad{} \biggr\}. 	  
	\end{equation}
	Therefore, \eqref{d-ramana} can be expressed as an SDP with auxiliary variables $W_i \in \mymatrix{n}$ 
	and $\beta_i \in \rad{}$ for $i=1, \dots, n+1;$ i.e., with polynomially many (in $n$ and $m$) variables and constraints. 
\end{Theorem}

We remark  that $\T(U)$ has a geometric meaning that is familiar from convex analysis; we will explain this 
in Section \ref{section-conclude}. However, to derive Ramana's dual, we only need the  algebraic definition of
$\T(U)$ in  \eqref{eqn-tan:2}.          

\co{
We note that $\T(U)$ has a geometric meaning,  that is familiar from convex analysis, as we will explain  
in Section \ref{section-conclude}. However, to derive Ramana's dual, we only need the  algebraic definition of
$\T(U)$ in  \eqref{eqn-tan:2}.          
}

\paragraph{Outline of the paper} 
Our proofs
use  mostly 
basic  linear algebra, but we  rely on a classic, rudimentary 
strong duality result from conic linear programming, that we recap in Proposition 
\ref{prop-Slater}. 

We employ the  notion of the 
 relative interior of a convex set. 
However, we  use it only  for a very particular type of set. 
For $k \geq 0$ we define 
\begin{equation} \label{eqn-C} 
\psd{k} \oplus 0  \, := \, \biggl\{  \bpx X & 0 \\ 0 &  0 \epx \in \sym{n} \, \Bigl|    X \in \psd{k} \,   \biggr\}. \, 
\end{equation}
It is straightforward to show that the relative interior of $\psd{k} \oplus 0$ is 
\begin{equation} \label{eqn-riC}
\biggl\{  \bpx X & 0 \\ 0 &  0 \epx  \in \sym{n} \, \Bigl|  \,  X \, \text{is } k \times k \, \text{positive  definite}    \biggr\},  
\end{equation}
and using  this fact the reader will be able to follow all proofs.  

In Section \ref{section-prelim} we introduce the key   ingredients of our proof:
i) {\em rescaling} the operator $\A$ and the right hand side $B; $ 
ii)  a {\em maximum rank  slack} in 
\eqref{p}, in other words, a maximum rank psd matrix of the form $B  - \A x, \,$where $x \in \rad{m};$
and iii)  Proposition  \ref{prop-Slater}.

Section \ref{section-proofs} has the main proofs, that we outline below: 
\begin{itemize}
	\item In Subsection 
\ref{subsection-gordan}, in Lemma  \ref{lemma-complete-certificate} (with  Lemma \ref{lemma-gordanstiemke} as a preliminary result)
 we describe certificates for  a maximum rank slack in \eqref{p}. 
By ``certificate" we mean 
a finite sequence of matrices that convince a third party that a maximum rank slack in \eqref{p} indeed has maximum rank. 

\item 
 In Subsection \ref{subsection-strongdual} we present a semidefinite program \eqref{rel-d}, which is a strong dual of \eqref{p}. 
 That is, \eqref{rel-d} is an ``inf" problem which attains its optimal value, when  finite, and 
 \begin{equation} \label{eqn-valP=valRelD} 
 	\val \eqref{p} = \val \eqref{rel-d}.
 \end{equation}
The problem \eqref{rel-d} is essentially the same as \eqref{d}, but in \eqref{rel-d}  only a block of the variable matrix $Y$ must be psd. 

At the same time, \eqref{rel-d} has a drawback: to write it down, we need to know a maximum rank slack  in \eqref{p} explicitly. 
However, in general we do not know such a slack explicitly; we only know that one exists. 


 \item 
 In Subsection    \ref{subsection-proof-of-theorem1}        we tie together the previously proved results and prove that 
  \eqref{rel-d} is ``mimicked" by Ramana's dual. We prove
\begin{equation} \label{eqn-val-Rel-D=val-Dram} 
	\val \eqref{rel-d} = \val \eqref{d-ramana},
\end{equation}
and that \eqref{d-ramana} attains its optimal value, when it is finite. 

In particular, to prove the inequality $\geq $ in \eqref{eqn-val-Rel-D=val-Dram} we produce 
an optimal solution of \eqref{d-ramana} as follows: from an optimal solution of \eqref{rel-d} we produce 
$U_{n+1}$ and $V_{n+1}, \, $ and from the certificates for a maximum rank slack given in Subsection \ref{subsection-gordan} we produce the other $U_i$ and $V_i.$ 

	We then  combine  \eqref{eqn-valP=valRelD} and  \eqref{eqn-val-Rel-D=val-Dram} and attainment in 
	\eqref{d-ramana} and prove Theorem \ref{theo:main}.

\item In Subsection \ref{subsection-thm2} we prove the SDP representation result Theorem \ref{theorem-repr}. 
	
\end{itemize}


In Section \ref{section-conclude} we conclude: we present a larger example and explain  the geometry of the set $\T(U), \, $ which at first may look enigmatic.

\section{Preliminaries}   
\label{section-prelim} 

\paragraph{Principal submatrices and concatenation}
Suppose $r$ and $s$ are integers in $\{1, \dots, n \},  r \leq s, \, $ and $Y \in \sym{n}.$ We then 
denote by $Y(r:s)$ the  principal submatrix  of $Y $ indexed by rows and columns $r, r+1, \dots, s.$ 

Further, we denote the  concatenation of  matrices $A$ and $B$  along the diagonal by $A \oplus B, \, $ i.e., 
$$
A \oplus B := \bpx A & 0 \\
0  &  B \epx.
$$

We naturally define the matrix $A \oplus B \oplus C$  as  $(A \oplus B) \oplus C.$ 

\subsection{Rescaling $\A$ and $B, \,$ slacks, and maximum rank slacks}

We will often rescale the operator $\A$ and the right hand side $B$ to put our semidefinite programs into a more convenient form. 
The precise definition follows:

\begin{Definition}
	We say that we {\em rescale} the operator $\A$ and the matrix $B$ if we perform the operations
	\begin{equation} \label{eqn-rotate} \tag{\mbox{Rescale}} 
		\ba{rcl}
		A_i & := & T^\top A_i T \, \text{for } \, i=1, \dots, m, \\ 
		B    & := & T^\top B T.
		\ena 
	\end{equation}
where $T$ is a suitable invertible matrix.
	\end{Definition}

\paragraph{Slacks and maximum rank slacks } 
We first define slack matrices in \eqref{p}, which generalize  slack vectors in linear programming.  
\begin{Definition}
	We say that $S \succeq 0$ is a {\em slack matrix} or {\em slack} in \eqref{p} if $S = B - \A x $ for some  
	$x \in \rad{m}.$ 
\end{Definition}
Note that 
if the $A_i$ and $B$ are diagonal, 
then \eqref{p} is a linear program, and the diagonal of a slack matrix  in \eqref{p} is just a slack vector in this linear program. 

Since the rank of a slack matrix is a nonnegative integer, and it is at most $n, \, $ 
the semidefinite program  \eqref{p} has a slack of maximum rank.
A maximum rank slack of \eqref{p}  will be a key player in the rest of the paper: it measures ``how centrally" the 
affine subspace $\{ \, B - \A x : x \in \rad{m} \, \} $ intersects the set of positive semidefinite matrices.

For convenience we  make the following assumption. 
\begin{Assumption}\label{assumption-slack} 
	There is a  maximum rank slack in \eqref{p} of the form 
	\vspace{.2cm} 
	$$
	Z \, = \, 
	\begin{pmatrix} I_r & 0 \\ 0 & 0 \end{pmatrix}, 
	$$
	where $0 \leq r \leq n,$ and, as usual, $I_r$ stands for an identity matrix of order $r.  \, $ 
\end{Assumption}
For the rest of this paper we fix this $Z$ and $r$.  

We can ensure Assumption \ref{assumption-slack} by a suitable rescaling as follows. Suppose  
$S$  is a maximum rank slack in \eqref{p},
$T$ is an invertible  matrix of suitably scaled eigenvectors of $S, \,$ and we 
perform the operations \eqref{eqn-rotate}. Afterwards 
in the new  primal problem a 
maximum rank slack is  $Z := T^\top S T, \, $ which is in the required shape. 

The maximum rank slack in \eqref{p} may not be unique. However, after we state 
Lemma \ref{lemma-complete-certificate}, we will prove a slightly weaker statement: any 
maximum rank slack in \eqref{p} must be of the form $R \oplus 0, \, $ where $R$ is an order $r$ symmetric positive definite matrix. 

\co{
We will prove this fact 
right after stating Lemma \ref{lemma-complete-certificate}. 
}

\co{ it is ``essentially" unique in the following sense: 
any 
maximum rank slack in \eqref{p} is of the form $R \oplus 0, \, $ where $R$ is an order $r$ symmetric positive definite matrix. 
We will prove this fact 
right after stating Lemma \ref{lemma-complete-certificate}. 
}

After the  initial rescaling that created our maximum rank slack $Z, \,$ we may rescale $\A$ and $B$ several more times. In these subsequent rescalings the $T$ transformation matrix will always be of the form 
	$$
T \, = \, 
\begin{pmatrix} I_r & 0 \\ 0 & M \end{pmatrix},
$$
where $M$ is an invertible, order $n-r$ matrix. These subsequent rescalings will keep $Z$ 
in the same form, since 
for any such $T$ matrix we have $T^\top Z T = Z. \, $ 

\begin{Example} (Example \ref{ex-main} continued)    
	As we discussed, in the SDP  \eqref{ex1-problem} the only feasible solution is $x_1=0, \, $ hence the maximum rank slack is just the right hand side
	\begin{equation} \label{eqn-Z-continued} 
		Z \, = \, \bpx 1 & 0 \\ 0 & 0 \epx.
	\end{equation}
That is, this SDP needs no rescaling, since  $Z$  already satisfies Assumption \ref{assumption-slack} (with $r=1$). 


\end{Example}

\subsection{Rescaling $\A$ and $B$ keeps the optimal value of Ramana's dual the same} 
\label{subsection-rescale-keeps-val-RamD-the-same}

When we rescale $\A$ and $B, \, $ 
we keep the optimal value 
(and  even the feasible set)   of \eqref{p} the same. It is a bit less obvious that  operations \eqref{eqn-rotate} do not affect 
the  optimal value  Ramana's dual,  so we prove that in Proposition \ref{prop-rotate-ramana}. 

First,  in Proposition \ref{prop-basic} we collect some useful properties of symmetric matrices. 

\begin{Proposition}
	\label{prop-basic}
	The following statements hold: 
	\begin{enumerate}
		\item \label{eqn-XYT} Suppose $X, Y \in \sym{n}$ and $T$ is  invertible. Then 
		\begin{equation} \nonumber   
			\la X, Y \ra \, = \, \la T^\top XT, T^{-1}YT^{-\top} \ra. 
		\end{equation} 
		\item 	\label{eqn-jerry} Suppose $U \in \psd{n}$ and $T$ is invertible. Then 
		\begin{equation}    \label{eqn-TtUT} 
			T^\top \T(U) T \, = \, \T(T^\top U T), 
		\end{equation} 
		where $T^\top \T(U) T$ is defined as $\{ T^\top V T \, | \, V \in \T(U) \,   \}.$ 
			\item \label{tan-shape} 
		Suppose $U \in \psd{n}$ is of the form $U = 0\oplus I_s, $ where  $s \leq  n$. 
		Then 
		$\T(U)$ is the set of matrices in $\sym{n}$ of the form
		
		\begin{equation} \label{eqn-tan-shape}
		\begin{pmatrix}[c|c]
			\bovermat{$n-s$}{{\mbox{$\,\, 0 \,\,$}}}	& \bovermat{$ s$}{{\mbox{$\,\,\,\ti \,\,\,\,\,$}}}  	\\ \hline 
			\ti & \ti 
		\end{pmatrix},
			\end{equation}
		
		where the elements in the $\ti$ blocks are arbitrary.
	\end{enumerate} 
\end{Proposition}  
\pf{} Statement \eqref{eqn-XYT} follows directly from the properties of the trace:
\beq \nonumber 
\ba{rcl}
\la X, Y \ra &  = &  \trace ( X Y  ) \\
&  = &  \trace ( X T T^{-1} Y T^{- \top} T^{\top} ) \\
& = &  \trace ( T^{\top} X T T^{-1} Y T^{- \top}) \\
& = & \la T^{\top} X T,  T^{-1} Y T^{- \top}  \ra. 
\ena 
\eeq
To prove   the inclusion $\subseteq$ in \eqref{eqn-jerry},  suppose $\R(W) \subseteq \R(U), $ so $W + W^\top \in \T(U). \, $ 
Then  
$$
\co{\R(W) \subseteq \R(U) \, \Rightarrow } \R(T^\top W T) \subseteq \R(T^\top U  T), 
$$ 
hence $T^\top ( W + W^\top) T \in \T(T^\top U T), \,$   so the inclusion $\subseteq$ follows. 

Next we apply the inclusion $\subseteq \, $ in  \eqref{eqn-TtUT},   but  now   with   $T^\top U T$ in place of $U$ and $T^{-1}$ in place of $T,$ 
and we deduce 
\begin{equation} \label{eqn-reverse}
	T^{-\top} \T(T^\top U T) T^{-1} \subseteq \T(U),
\end{equation} 
and left multiplying \eqref{eqn-reverse} by $T^\top$ and right multiplying  by $T$ yields the inclusion $\supseteq$ in \eqref{eqn-jerry}.

\myblue{
The statement \eqref{tan-shape} follows  from  the definition of the tangent space in \eqref{eqn-tan:2} and since for 
 a matrix $W \in \rad{n \times n}$ we have $\R(W) \subseteq \R(U)$ exactly when  the first 
$n-s$ rows of $W$ are zero.}
\qed

To build intuition, we first argue that rescaling keeps the optimal value of \eqref{d} the same. 
Indeed, suppose $Y$ is feasible in \eqref{d} with objective value, say, $\alpha, \, $ 
and we apply the operations \eqref{eqn-rotate} using an invertible matrix $T.$
Then  $T^{-1} Y T^{- \top}$ is feasible in \eqref{d} after rescaling, and has objective value  $\alpha: \!$ this follows  by item \eqref{eqn-XYT} in Proposition \ref{prop-basic}, with $A_i$ or $B$ in place of $X.  \,$ 
Since we can undo the rescaling with $T$ by another rescaling (with $T^{-1} \, ), \,$ it follows that the optimal value of \eqref{d} 
is the same before and after rescaling.

A similar argument, given in Proposition  \ref{prop-rotate-ramana}, shows that the optimal value of \eqref{d-ramana} stays the same after rescaling. The  only difference is that we now also have to take care of the tangent space constraints \eqref{eqn-tan}. 

\begin{Proposition} \label{prop-rotate-ramana}  
	The operations \eqref{eqn-rotate} keep the optimal value of \eqref{d-ramana} the same.
\end{Proposition}
\pf{} Suppose that $(U_j, V_j)_{j=0}^{n+1}$ is feasible in \eqref{d-ramana} with objective value $\alpha$ before we performed 
\eqref{eqn-rotate}; we prove that  
$$(U_j', V_j')_{j=0}^{n+1} := (T^{-1} U_j T^{-\top}, T^{-1} V_j T^{-\top})_{j=0}^{n+1}$$ is feasible  afterwards, and has the same objective value. 

Indeed, by \eqref{eqn-XYT} in Proposition \ref{prop-basic} we have 
\beq \nonumber 
\ba{rcl}
\la A_i, U_j + V_j \ra & = & \la T^\top A_i T, U_j' + V_j' \ra
\ena
\eeq
for all $i$ and $j. \, $ 
 By the same logic we see that 
\beq \nonumber 
\ba{rcl}
\la B, U_j + V_j \ra & = & \la T^\top B T, U_j' + V_j' \ra 
\ena
\eeq
for all $j.$ 
Thus, $(U_j', V_j')_{i=0}^{n+1}$  satisfies the equality  constraints of \eqref{d-ramana}
after we executed \eqref{eqn-rotate}. 

Also, 
for all $j$  we have $V_j \in \T(U_{j-1}).$ Thus,  by item \eqref{eqn-jerry} in Proposition \ref{prop-basic} (with   $T^{-\top}$ in place of $T$) we deduce 
$$
V_j' \in \T(U_{j-1}') \, 
$$
for all $j.$ 
 
Summarizing,  $(U_j', V_j')_{j=0}^{n+1}$   is feasible in  \eqref{d-ramana} with objective value $\alpha$ 
after we executed \eqref{eqn-rotate}, completing the proof. \qed 

\subsection{Strong  duality under Slater's condition} 
We next state a classic strong duality result assuming the underlying space is $\symn, \, $ a special 
case most relevant for this work \footnote{Proposition \ref{prop-Slater} is usually stated in the space 
	$\rad{n}.$}. 
 
Suppose $K \subseteq \symn$ is a closed convex cone \footnote{That is, $K$ is closed, convex, and $\lambda x \in K$ for all $x \in K$ and $\lambda \geq 0.$} and let us 
denote its relative interior by $\ri K, \, $ and its dual cone by $K^*,$ i.e., 
\begin{equation} \nonumber
K^* \, = \, \{  Y \in \symn \, | \, \inProd{X}{Y}  \geq 0 \, \forall \, X \in K \,  \}.
\end{equation}
\bprop \label{prop-Slater} 
Suppose $\A, B, K,  \, $ and $K^*$ are as previously defined, and 
\begin{equation} \label{eqn-slater} 
B - \A x \in \ri K \, \text{for some } \, x \in \rad{m}. 
\end{equation}   
Then  
\begin{equation} \label{eqn-p-d-C} 
\sup \{ \, c^\top x  \mid B - \A x \in K  \} \, = \, \inf \{ \inProd{B}{Y} \mid \A^*Y = c, Y \in K^*\}
\end{equation} 
and the optimal value of the ``inf" problem is attained when it is finite. 
 \qed 
\eprop

When condition \eqref{eqn-slater} holds, we say that the ``sup" problem in 
\eqref{eqn-p-d-C} satisfies Slater's condition. 

For better intuition, we next outline  two important uses of Proposition \ref{prop-Slater}. Each one corresponds  to how large $r, \, $ the rank of $Z$ is. 

First suppose $r = n,  \, $ and we choose $K = \psd{n}.$ Then $K^* = \psd{n}$  holds as well, so the optimization problems in 
\eqref{eqn-p-d-C} are just the semidefinite programs  \eqref{p} and \eqref{d}. Further,  \eqref{p} satisfies Slater's condition.
Hence \eqref{d} has the same optimal value as \eqref{p}, and attains this value when it is finite. In other words, the usual dual is just as good as Ramana's dual.

Second, suppose $r < n$ and we set $K = \psd{r} \oplus 0.$ Since the dual cone of $\psd{r}$ is $\psd{r}, \, $ itself, we have 
\begin{equation} \label{eqn-dual-psdk} 
	K^* \, = \, \{   Y \in \sym{n} \, | \, Y(1:r) \succeq 0 \},  
\end{equation}
i.e., in the dual cone only the upper left $r \times r$ block of matrices must be psd. Suppose   in the constraint set of \eqref{p}, namely in 
$$
B - \A x \in \psd{n} 
$$
we replace $\psd{n}$ by $K$ and in the constraint set of \eqref{d} we replace $Y \in \psd{n}$ by $Y \in K^*.$ 
Then we show in Lemma \ref{lemma-relaxed-dual} that three interesting things happen.
First, even though $K$ is a smaller set than $\psd{n}, \, $ the feasible set of 
\eqref{p} remains the same. Second, \eqref{p} satisfies Slater's condition. Third, by Proposition 
\ref{prop-Slater}, the dual becomes  the promised strong dual \eqref{rel-d}.  

A bit surprisingly, short and self-contained proofs of Proposition~\ref{prop-Slater} 
are rare in the literature.  Fortunately, such a 
proof is given in Theorem 7 in the technical 
report by Luo, Sturm and Zhang \cite{LuoSturmZhang:97}. 
However,  proofs of more specific or more general 
statements are common. 
As to more specific ones, the result with 
``relative interior" replaced by ``interior" appears in 
Section 3.2 of Renegar \cite{Ren:01},   in Section 2.4 of the textbook of Ben-Tal and Nemirovskii \cite{BentalNem:01}, 
and in Section 5.3 of  Borwein and Lewis \cite{BorLewis:05}. 
As to more general statements, Proposition~\ref{prop-Slater}  
follows from Fenchel's duality theorem 
 in Rockafellar's classic text \cite[Theorem~31.4]{Rockafellar:70}.

\section{Proofs} 
\label{section-proofs} 

	Recall that whenever \eqref{p} is feasible, it has a maximum rank slack. Further, after a suitable rescaling, we fixed a maximum rank slack $Z$ 
	(with  rank $r$) 
of the shape given in Assumption \ref{assumption-slack}. This $Z$ will play a key role in all proofs.

\subsection{Certificates for the maximum rank slack} 
\label{subsection-gordan} 

	Lemma \ref{lemma-gordanstiemke} below proves that all slacks
	of \eqref{p} have certain restrictions on their shape and rank.

\begin{Lemma} \label{lemma-gordanstiemke} 
	Suppose $s \in \{r+1, \dots, n \}$ is an integer.  
	Then the following semidefinite system has a solution: 
	\beq \label{eqn-Astar-Yj+1} 
	\ba{rcl}
	\A^* Y & = & 0 \\
	\la B, Y \ra & = & 0 \\
	Y & \in & \sym{n} \\
	Y(1:s) & \in & \psd{s} \setminus \{ 0 \}.
	\ena
	\eeq 
\end{Lemma} \qed

The $Y$ matrix given in Lemma \ref{lemma-gordanstiemke} certifies that 
	\eqref{p} has no slack whose order $s$ leading 
	principal submatrix  is positive definite, and the rest is zero. Indeed, suppose  $S = B - \A u$ is such a slack, then we 
have 
\begin{equation} \label{eqn-SY=0} 
	0 \,=  \,  \la B, Y \ra - \la \A^* Y, u \ra  \, = \, \la B - \A u, Y  \ra \, = \, \la S, Y \ra. 
\end{equation}
However, $\la S, Y \ra > 0, \, $ which is a contradiction. 
\bex (Example \ref{ex-main} continued) 
As we previously discussed, in the SDP \eqref{ex1-problem} the right hand side 
$ Z \, = \, 
\bigl( \begin{smallmatrix} 1         & 0\\
	0    & 0 \end{smallmatrix} \bigr)
$
is the maximum rank slack, 
so $r=1.$  
By Lemma \ref{lemma-gordanstiemke} with $s=2$ we produce  
\begin{equation} \label{eqn-Y-cert} 
	Y = \bpx 0 & 0 \\ 0 &  1 \epx
\end{equation}
to certify  that \eqref{ex1-problem} has no positive definite slack. 
In  this tiny example, 
since $Y$ in \eqref{eqn-Y-cert} certifies that \eqref{ex1-problem} has no rank $2$ slack, it  also certifies  that $Z$ has maximum rank.

However, in larger examples  we will need a finite {\em sequence}  of matrices to completely certify that $Z$ has maximum rank. We will illustrate this fact in  Example \ref{ex-larger}.
\eex

\pf{of Lemma \ref{lemma-gordanstiemke}} 
Let $K = \psd{s} \oplus 0.$ Consider the primal-dual pair of semidefinite programs
\begin{equation} \label{eqn-Gordan} 
	\sup_{x,t}  \{ \, t   \mid B - \A x - t ( I_s \oplus 0) \in K  \}  \; \text{and} \; \inf_Y  \{ \la B, Y \ra \mid \A^*Y = 0, \la ( I_s \oplus 0), Y \ra = 1, \, Y \in K^*      \},
\end{equation} 
and for brevity, define 
$$ 
S(x,t) := B - \A x -t  ( I_s \oplus 0 ) \; \text{for} \; x \in \rad{m} \; \text{and} \; t \in \rad{}.
$$
We first claim that the optimal values of the ``sup" and ``inf" problems in \eqref{eqn-Gordan} are the same and the optimal value of the ``inf" problem is attained when it is finite.
For that, let $x \in \rad{m}$ be such that $Z = B - \A x. $ Then the upper left $s \times s$ block of 
$S(x,-1)$ is positive definite, and the other elements of $S(x,-1)$ are zero. Thus $S(x, -1) \in \ri K, \, $ so the 
``sup" problem in \eqref{eqn-Gordan} satisfies Slater's condition, hence 
our claim follows from Proposition \ref{prop-Slater}.

Next we claim that the optimal value of both optimization problems in \eqref{eqn-Gordan}  is nonnegative. For that, again let 
$x \in \rad{m}$ be such that $Z = B - \A x, $  then $S(x,0) \in K, $ so  the optimal value of the ``sup" problem is indeed nonnegative, and 
our claim follows. 

Then we claim that the optimal value of both optimization problems in \eqref{eqn-Gordan}  is zero. 
To obtain a contradiction, suppose that $S(x,t) \in K $ for some $x \in \rad{m}$ and $t > 0.$ Then the upper left $s \times s$ block of 
$B - \A x$ is positive definite, and the other elements are zero. 
Thus $B - \A x$ is a slack in \eqref{p} whose rank is larger than $r, \, $ which is the required contradiction.

Thus the ``inf" problem in \eqref{eqn-Gordan} has a feasible solution $Y$ with objective value zero.
By 
$Y \in K^*$ we get $Y(1:s) \succeq 0$ (see \eqref{eqn-dual-psdk}) and by  
$ \la \bigl(  I_{s}  \oplus 0 \bigr), Y \ra =1$ we get $Y(1:s) \neq 0.$ Thus $Y$ 
satisfies \eqref{eqn-Astar-Yj+1}, as wanted. 
\qed

Lemma \ref{lemma-gordanstiemke} gave  a partial certificate for the maximum rank slack $Z:$ the $Y$ matrix in 
\eqref{eqn-Y-cert} certifies that $Z$ has maximum rank among slacks with a fixed form. 
In contrast, Lemma \ref{lemma-complete-certificate} gives a complete certificate: it shows that $Z$ has maximum rank among slacks of any form. 

\begin{Lemma} \label{lemma-complete-certificate} 
	We can rescale $\A$ and $B$ so that after the rescaling the following hold:
	\begin{enumerate}
		\item \label{lemma-complete-certificate-1} The 
		$Z$ matrix given in Assumption \ref{assumption-slack} is still a maximum rank slack in \eqref{p}. 
		\item There exist symmetric matrices $Y_1, \dots, Y_k$ which are of the form
	
	\vspace{.2cm} 
		\beq \label{eqn-Yi} 
		\ba{rclrcl}
		Y_i \!\!\!\! & := & \!\!\!\! 
		\begin{pmatrix}[c|c|c]
			\bovermat{$n - \sum_{\ell=1}^{i} r_\ell$}{\phantom{00000000000}} & \bovermat{$r_i$}{\phantom{0000}}	& \bovermat{$\sum_{\ell=1}^{i-1} r_\ell$}{\,\,\,\,\,\,\,\,\,\,\,\,\, \times \,\,\,\,\,\,\,\,\,\,\,\,\,\,}	\\ \hline 
			\pha{0}  & I  & \times \\ \hline 
			\times  & \times & \times 
		\end{pmatrix}, \!\!\!   
	\end{array} 
\eeq 
and satisfy  
\beq \label{eqn-Yi-A*} 
\ba{rcl} 
\A^{*}  Y_i  & = &  0 \\
\la B, Y_i \ra   & = & 0
\ena
\eeq 
for $i=1, \dots, k.$ 
Here $k \geq 0$ and the $r_i$ are positive integers such that 
$\sum_{i=1}^k r_i = n-r.$ 
	\end{enumerate}
\qed	
\end{Lemma}

Here, and in the sequel empty blocks in matrices contain all zeros, and  $\times$ blocks may have arbitrary elements.
(In some matrices we still explicitly indicate zero entries, if this helps readability.)  

How do the $Y_i$ in Lemma \ref{lemma-complete-certificate} certify that  $Z$ has maximum rank? 
To explain, let $S = B  - \A u$ be an arbitrary slack in \eqref{p}, where $u \in \rad{m}.$  
We present a simple argument using the $Y_i$ to prove 
that the rank of $S$ is at most $r, \, $ in particular, 
that the last $n-r$ rows and columns of $S$ are zero.

Using an argument similar to the one in  \eqref{eqn-SY=0} we first deduce $\la S, Y_1 \ra = 0.$
Since $\la S, Y_1 \ra$ is the sum of the last $r_1$ diagonal elements of $S, \, $ and these elements  are all nonnegative,  they must be all zero. 
Since $S \succeq 0, \, $ we learn that 
the last $r_1$ rows and columns of $S$ are zero. 

We then repeat the above argument with $Y_2$ in place of 
$Y_1. \, $ We have $\la S, Y_2 \ra = 0.  \, $ Since the last $r_1$ rows and columns of $S$ are zero, $\la S, Y_2 \ra$ is 
the sum of the diagonal elements of $S$ in rows numbered 
${n-r_1-r_2+1}, \dots, {n-r_1}.$ So these rows and colums are all zero. 

Continuing, since the sum of all $r_i$ is $n-r, \, $ 
we deduce that the last $n-r$ rows and columns of $S$ are zero, as required. 

This argument also proves that $Z$ is a maximum rank slack that is unique  up to rescaling. Precisely,  it proves that any rank $r$ slack in \eqref{p} 
must look like $R \oplus 0, \, $ where 
$R$ is order $r$ and positive definite.

Continuing Example \ref{ex-main}, the $Y$ matrix in  \eqref{eqn-Y-cert} can serve as $Y_1$ for the SDP \eqref{ex1-problem}:  
here we can choose $k=1$ so there is no need for other $Y_i.$

\pf{of Lemma \ref{lemma-complete-certificate}} 
For a nonnegative integer $j$ we 
consider the following conditions: 

\begin{enumerate} 
	
		\item \label{eqn-Z} 
		We have rescaled $\A$ and $B$ so that $Z$ is still a maximum rank slack in \eqref{p} after the rescaling.
		
	\item \label{eqn-construct-Yi}   
 We have   constructed $Y_i$ for $i=1, \dots, j$ 
which are of the form required in \eqref{eqn-Yi} 
and satisfy the equations \eqref{eqn-Yi-A*}. 
Further,  the $r_i$ sizes of the identity blocks in the $Y_i$ matrices are all positive.
\end{enumerate} 

We start with $j=0,$ then \eqref{eqn-Z} is satisfied by assumption, and \eqref{eqn-construct-Yi} is satisfied vacuously. 
\co{both invariant conditions are satisfied. 

We start with $j=0,$ then both invariant conditions are satisfied. 

}

In a general step  we assume that $j \geq 0 \,$ and that conditions \eqref{eqn-Z} and  \eqref{eqn-construct-Yi} 
hold with $j. \,$ The argument after the statement of Lemma 
\ref{lemma-complete-certificate} implies that the identity blocks in the $Y_i$ and the identity block in $Z$ do not overlap, hence 
$ \sum_{i=1}^j r_i \leq n-r.$  

We define
$$
s := n - \sum_{i=1}^j r_i, 
$$
hence $s \geq r.$ 

If $s = r, \, $ we let $k=j, \, $ and stop. 

If $ s >r,$ then we will construct matrix $Y_{j+1},$ make sure that conditions \eqref{eqn-Z} and  \eqref{eqn-construct-Yi} hold with $j+1 \, $ in place of $j, \, $ 
and increment $j.$ 
Since all  the $r_i$ are positive, we can execute this step at most $n$ times.

First we invoke Lemma \ref{lemma-gordanstiemke}  and produce $Y_{j+1} \in \symn$ such that 
\beq \nonumber 
\ba{rcl} 
\A^{*}  Y_{j+1}   & = &  0 \\
\la B, Y_{j+1}  \ra   & = & 0 \\
Y_{j+1}(1:s) & \in  & \psd{s} \setminus \{ 0 \}.  
\ena
\eeq 
Let us recall the maximum rank slack $Z$ from Assumption \ref{assumption-slack}. Then using an argument just like 
in \eqref{eqn-SY=0}, we deduce $\la Z, Y  \ra  \, = \,  0. $ 
Thus the first $r$ rows and  columns of $Y_{j+1}(1:s)$ are zero, meaning $Y_{j+1}$ looks like 

\vspace{.1cm} 		 			
\beq \nonumber 
Y_{j+1}   \, = \,   
\begin{pmatrix}[c|c|c]
	\bovermat{$r$}{\phantom{000}} & \bovermat{$s-r$}{\phantom{0000000}}	& \bovermat{$n-s$}{\pha{00}\ti\pha{000}} 	\\ \hline 
	\pha{0} & \bar{Y}   & \ti  \\ \hline 
	\ti  & \ti  & \ti 
\end{pmatrix}, \,\,  
\eeq 
\vspace{.1cm} 
where $\bar{Y} \succeq 0 \, $ and,  as usual,  the elements in the $\ti$ blocks are  arbitrary.

We next put  $Y_{j+1}$ into the required format, and rescale $\A$ and $B$ to make sure 
that conditions \eqref{eqn-Z} and  \eqref{eqn-construct-Yi} hold with  $j+1$ in place of $j.$  
For that, let $r_{j+1}$ be the rank of $\bar{Y},$  
$Q$  be an invertible  matrix of eigenvectors of 
$\bar{Y}$ such that $Q^\top \bar{Y} Q = 0 \oplus I_{r_{j+1}}$ and 
\begin{equation} \label{eqn-T-Q} 
	T := I_r \oplus Q \oplus I_{n-s}. 
\end{equation}
We next perform the operations
\beq \label{eqn-rotate-2} 
\ba{rcll} 
Y_i & := & T^\top Y_i T               &  \text{for} \,\, i=1, \dots, j+1. 
\ena
\eeq
These operations put $Y_{j+1}$ into the form required in \eqref{eqn-Yi} , and keep $Y_1, \dots, Y_j$ in the same form. 

However, now the $Y_i$ may not satisfy equations \eqref{eqn-Yi-A*}, so we correct that issue next. 
We perform the operations 
\beq \label{eqn-rotate-3}
\ba{rcll} 
A_i & := & T^{-1} A_i T^{-\top}  & \text{for} \,\, i =1, \dots, m,  \\ 
B & := & T^{-1} B T^{-\top},
\ena 
\eeq
i.e., the \eqref{eqn-rotate} operations with $T^{- \top}$ in place of $T.$ After these operations the equations \eqref{eqn-Yi-A*} hold  (by part \eqref{eqn-XYT} of Proposition \ref{prop-basic}), 
so condition \eqref{eqn-construct-Yi}  holds with $j+1.$ 

Finally we show that condition 
\eqref{eqn-Z} remains true. For that, we  observe  that after performing the operations 
in \eqref{eqn-rotate-3}, the matrix $T^{-1} Z T^{- \top} $ is a slack of rank $r$ in \eqref{p}, hence it is a maximum rank slack.
Also, by \eqref{eqn-T-Q} we have 
\begin{equation} \label{eqn-T-Q-inverse} 
	T^{-1} := I_r \oplus Q^{-1} \oplus I_{n-s},
\end{equation}
hence $T^{-1} Z T^{- \top} = Z, \, $ so $Z$ is a maximum rank slack in \eqref{p}. Thus \eqref{eqn-Z} holds, and the proof is complete.
\qed

The proof of Lemma \ref{lemma-complete-certificate} gives a {\em facial reduction algorithm} to construct the $Y_i.$ 
To explain this parlance, note that the set 
$$
F := \psd{r} \oplus 0
$$
is a {\em face} of $\psd{n} \, $ \footnote{This means two things: (i)  it is a convex subset of $\psd{n}$ and (ii) if $X, Y$ are in $\psd{n},$  and 
	the open line segment $\{ \lambda X + (1-\lambda)Y : \lambda \in (0,1) \}$ intersects $F, \, $ then both $X$ and $Y$ must be in $F.$} and the $Y_i$ matrices {\em reduce} the set of feasible slacks of \eqref{p} 
to live in $F. \, $ 
For simplicity, in our proofs 
we do not mention faces; there is no need, since $F$ is perfectly captured by the maximum rank slack $Z.$ 

\myblue{
We note that the algorithm is theoretical, since to implement it, we must find the $Y_i$ certificates in Lemma \ref{lemma-complete-certificate} in exact arithmetic;
and for that, we must solve the pair of semidefinite programs in \eqref{eqn-Gordan} in exact arithmetic. However, some  heuristic
 implementations of facial reduction 
exist,  see for example \cite{permenter2014partial} and \cite{zhu2019sieve}. }
We further refer the reader to \cite{BorWolk:81, Pataki:13, WakiMura:12} 
for facial reduction algorithms for more general problems.

\subsection{A strong dual, assuming we know a maximum rank slack}
\label{subsection-strongdual} 

In this subsection we present our promised strong dual \eqref{rel-d}.
\begin{Lemma} \label{lemma-relaxed-dual} 
	Suppose that we rescaled $\A$ and $B$ as stated in Lemma \ref{lemma-complete-certificate}, and consider the optimization problem 
	\begin{equation}\label{rel-d}
		\begin{split}
			\inf  & \,\, \la B,  Y  \ra  \\
			s.t. & \,\, \A^* Y = c  \\ 
			& \,\,   Y \in (\psd{r} \oplus 0)^*. 
		\end{split}\tag{\mbox{${\rm D_{\rm strong}}$}}
	\end{equation}
Then  
$$
\val \eqref{p} = \val \eqref{rel-d},
$$
and the optimal value of \eqref{rel-d} is attained when it is finite. 
\end{Lemma}
\pf{} By the argument after the statement of Lemma \ref{lemma-complete-certificate} we see that any 
slack in \eqref{p} is contained in $\psd{r} \oplus 0.$ Thus \eqref{p} is equivalent to 
\begin{equation} \label{rel-P} 
	\sup \{ \, c^\top x  \mid B - \A x \in \psd{r} \oplus 0 \, \}.
\end{equation}
Again recall the maximum rank slack $Z$ from Assumption \ref{assumption-slack}. 
Since $Z$ is in the relative interior of $\psd{r} \oplus 0, \, $ we see that 
\eqref{rel-P} satisfies Slater's condition. Thus by Proposition \ref{prop-Slater} the dual of \eqref{rel-P} attains its optimal value, when it is finite,
and this optimal value is the same as the optimal value of \eqref{rel-P}. 

But the dual of \eqref{rel-P} is just \eqref{rel-d}, hence the proof is complete.
\qed

We next illustrate Lemma \ref{lemma-relaxed-dual}.
\begin{Example} (Example \ref{ex-main} continued)
	We repeat the SDP from Example \ref{ex-main} for convenience.
		\beq \label{ex1-problem-2}
	\begin{array}{rl}
		\sup &   2 x_1 \\
		s.t. & x_1 \bpx 0 & 1 \\ 1 & 0 \epx \preceq \bpx 1 & 0 \\ 0 & 0 \epx.
	\end{array} 
	\eeq
	As we discussed, the only feasible solution is $x_1=0, \, $ hence the optimal value of 
	\eqref{ex1-problem-2} is zero, and the right hand side in  \eqref{ex1-problem-2} is the maximum rank slack.

	Thus in the strong dual \eqref{rel-d} of \eqref{ex1-problem-2} only the upper left $1 \times 1$ block of $Y$ must be psd. Hence 
	\eqref{rel-d} is 
	 	\beq \label{ex1-problem-dual-strong} 
	\ba{rrcl}
	\inf & y_{11} \\
	s.t. & y_{12}  & = &  1 \\
	& y_{11} & \geq & 0,
	\ena
	\eeq
	which is  just a linear program.
	The matrix 
	$$
	Y^* = \bpx 0 & 1 \\ 1 & 0 \epx
	$$
	is an optimal solution of \eqref{ex1-problem-dual-strong} that attains the optimal value of $0.$ 
	\end{Example} 

We note that the strong dual \eqref{rel-d} is essentially the same as the usual dual \eqref{d}, however, 
it requires only a block of the variable matrix $Y$ to be psd. 
Thus \eqref{rel-d} should even be easier to solve than \eqref{d}!  So why not use it?     

Here is the catch: to write down \eqref{rel-d} we would need to know a maximum rank slack in \eqref{p} explicitly.
If we did, then by rescaling we could ensure 
that a maximum rank slack $Z \, $ is in the shape required in Assumption 
\ref{assumption-slack}, then write down \eqref{rel-d}. Of course, in general we do not know a maximum rank slack explicitly, we only know 
that one exists.  

However, in the next section we show that \eqref{rel-d} is ``mimicked"  
by Ramana's dual, which has no need of a maximum rank slack;  of course, Ramana's dual needs many more variables. 

\co{
we do know $Z, \, $  implicitly, using the $Y_i$ certificates from Lemma \ref{lemma-complete-certificate}.
In Lemma \ref{lemma-valP>=valDram} we will use these $Y_i$ to construct a feasible solution of Ramana's dual with objective value at least as good 
as that of \eqref{rel-d}. 
}

\subsection{Proof of Theorem \ref{theo:main}} 
\label{subsection-proof-of-theorem1}

In this subsection we complete the proof of Theorem \ref{theo:main}. 
First,  in Lemmas \ref{lemma-valP>=valDram}  and  \ref{lemma-valP<=valDram} we prove 
\begin{equation} \label{eqn-kramer} 
	\val \eqref{rel-d} = \val \eqref{d-ramana},
\end{equation}
and if $\val \eqref{rel-d}$ is finite, then \eqref{d-ramana} has a solution with that value. 

\begin{Lemma} \label{lemma-valP>=valDram} 
	$$
	\val \eqref{rel-d}  \geq \val \eqref{d-ramana}. 
	$$
Further, when $\val \eqref{rel-d}$ is finite, \eqref{d-ramana} has a solution with that value.
\end{Lemma}

\pf{} \myblue{If \eqref{rel-d} is infeasible, then there is nothing to prove, so let us assume it is feasible. }
Further, assume that we rescaled $\A$ and $B$ as stated in Lemma \ref{lemma-complete-certificate}. 
By Lemma \ref{lemma-relaxed-dual} we have that \eqref{rel-d} has an optimal solution, and we choose 
$Y^* \in (\psd{r} \oplus 0)^*$  to be an optimal solution.

Let  $Y_1, \dots, Y_k$ be the matrices we constructed in Lemma \ref{lemma-complete-certificate}. Recall that $k \leq n.$ 
We will construct a feasible solution of \eqref{d-ramana} with value equal to $\val \eqref{rel-d}.$ 

First we outline the idea. From \eqref{eqn-Yi} we see that  each $Y_i$ can be written as $Y_i = U_i + V_i, \, $ where $U_i$ is psd, and 
$V_i \in \tan(U_{i-1}).$ Precisely,  we may choose $U_i := 0 \oplus I_{r_1 + \dots + r_i}$  for all $i,$ 
then $V_i \in \tan(U_{i-1})$ follows from part \eqref{tan-shape} of Proposition \ref{prop-basic}.

Then 
we can decompose $Y^*$ as 
$$Y^* \, = \, U_{k+1} + V_{k+1}, \, \text{with} \, U_{k+1} \in \psd{n}, \, \text{ and} \, V_{k+1} \in \tan(U_k). $$ 
In particular, we can choose $U_{k+1} := Y^*(1:r) \oplus 0$ 
and $V_{k+1} := Y^* - U_{k+1}.$ Then  $V_{k+1}  \in \tan(U_{k})$ follows from part \eqref{tan-shape} of Proposition \ref{prop-basic} and 
$r_1 + \dots + r_k = n - r.$ 

Thus, if we set $U_0 = V_0 = 0, \, $  we obtain a feasible solution to a variant of \eqref{d-ramana} in which $n$ is replaced by $k.$

\co{
\myblue{
First we outline the idea. 
Suppose we let $U_i := 0 \oplus I_{r_1 + \dots + r_i}$ and $V_i := Y_i - U_i$ for  $i =1, \dots, k.$ 
Then  $Y_i = U_i + V_i,$ and from \eqref{eqn-Yi}
we have $V_i \in \tan(U_{i-1}).$ 
In addition, suppose we let $U_{k+1} = Y^*(1:r) \oplus 0$ 
and $V_{k+1} = Y^* - U_{k+1}.$ Since  $U_k = 0 \oplus I_{n-r}$,
we have 
\[
Y^* \, = \, U_{k+1} + V_{k+1}, \, \text{with} \, U_{k+1} \in \psd{n}, \, \text{ and} \, V_{k+1} \in \tan(U_k). 
\]
Thus we obtain a feasible solution to a variant of \eqref{d-ramana} in which $n$ is replaced by $k.$ 
}
}

This plan is not quite perfect, since 
$k$ may be strictly less than $n. \, $  If this happens,  we need to modify our plan, namely we need 
to add some zero $U_i$ and $V_i$ at the start to create a feasible solution of \eqref{d-ramana}.

We now carry  out this modified plan. 

\begin{enumerate}
	\item The first few $U_i$ and $V_i$ are ``padding": we set 
	$$
	U_0 = V_0 = \dots = U_{n-k} = V_{n-k} = 0.
	$$
	Then 
\eqref{eqn-ai0}, \eqref{eqn-b0},  \eqref{eqn-tan}, and   \eqref{eqn-Uipsd} hold in Ramana's dual for $i \leq n-k.$

\item Then from $Y_1, \dots Y_k$  we construct $U_{n-k+1}, V_{n-k+1}, \dots, U_n, V_n.$ 
We write 

\vspace{.1cm}
\beq \nonumber 
Y_i \, = \,  
\underbrace{\begin{pmatrix}[c|c|c]
	\bovermat{$n - \sum_{\ell=1}^{i} r_\ell$}{\phantom{0000000000}} & \bovermat{$r_i$}{\phantom{0000000000}}	& \bovermat{$\sum_{\ell=1}^{i-1} r_\ell$}{\phantom{00000000000}}	\\ \hline 
	& I  & \\ \hline 
	\phantom{00000000000}   & \phantom{00000000000} & I 
\end{pmatrix}}_{U_{n-k+i}} + \underbrace{\begin{pmatrix}[c|c|c]
	\bovermat{$n - \sum_{\ell=1}^{i} r_\ell$}{\phantom{0000000000}} & \bovermat{$r_i$}{\phantom{0000000000}}	& \bovermat{$\sum_{\ell=1}^{i-1} r_\ell$}{\,\,\,\,\,\,\,\,\,\,\,\,\,\, \times \,\,\,\,\,\,\,\,\,\,\,\,\,\,\,}	\\ \hline 
	\pha{0}  &  & \times \\ \hline 
	\times  & \times & \times  
\end{pmatrix}}_{V_{n-k+i}} 
\eeq 
for $i=1, \dots, k.$ In other words, we let $U_{n-k+i}$ as above, then set 
$V_{n-k+i} := Y_i - U_{n-k+i}. $   
\myblue{ Then  by Part  \eqref{tan-shape}  in Proposition \ref{prop-basic} we have 
$$V_{n-k+i} \in \tan(U_{n-k+i-1}) \; \text{for} \; i=1, \dots, k,$$ } 
so
	\eqref{eqn-ai0}, \eqref{eqn-b0}, \eqref{eqn-tan}, and \eqref{eqn-Uipsd}   hold  in Ramana's dual  for $i=n-k+1, \dots, n.$ 
(It is useful to note that $V_{n-k+1}=0.$ )


\item We finally  split $Y^*$ to construct $U_{n+1}$ and $V_{n+1}:$ 
we write $Y^*$ as 
 
  \vspace{.1cm} 
 \beq \nonumber 
 Y^* = \underbrace{\bpx Y^*(1:r) & 0 \\ 0 & 0 \epx}_{U_{n+1}} + 
 \underbrace{\begin{pmatrix}[c|c]
 	\bovermat{$r$}{{\mbox{$\,\, 0 \,\,$}}}	& \bovermat{$n - r$}{{\mbox{$\,\,\,\ti \,\,\,\,\,$}}}  	\\ \hline 
 	\ti & \ti 
 \end{pmatrix}}_{V_{n+1}}.
 \eeq
 That is, we let $U_{n+1}$ be as above, then set $V_{n+1} := Y^* - U_{n+1}.$ 
 Then $\A^* (U_{n+1} + V_{n+1}) = c \,$  and $U_{n+1} \in \psd{n}.$ Also, by Part \eqref{tan-shape} of Proposition \ref{prop-basic} and $U_n = 0 \oplus I_{n-r}$ 
 we have $V_{n+1} \in \tan (U_n).$ 

\end{enumerate}  
In summary, 
$U_0, V_0, \dots, U_{n+1}, V_{n+1}$ 
is a feasible solution of \eqref{d-ramana} with value equal to 
$\val \eqref{rel-d}, \, $ hence the proof is complete.
\qed

We invite the reader to follow the recipe in the proof above, and construct the optimal solution of the Ramana dual of 
the problem \eqref{ex1-problem}. This solution was already given in Example \ref{ex1-continued}, but it is fruitful to produce it from 
the following ingredients: the maximum rank slack $Z$ which is just the right hand side in    \eqref{ex1-problem};
and the  $Y = Y_1 $ matrix in  in  \eqref{eqn-Y-cert}  that certifies that $Z$ indeed has maximum rank. 

We next prove the inequality $\leq$ in \eqref{eqn-val-Rel-D=val-Dram}. 

\begin{Lemma} \label{lemma-valP<=valDram} 
	We have 
	\begin{equation} \label{elaine} 
	\val \eqref{rel-d}  \leq \val \eqref{d-ramana}.
	\end{equation}
\end{Lemma}
\pf{} \myblue{ If \eqref{d-ramana} is infeasible, then there is nothing to prove, so assume it is feasible, }
and let $U_0, V_0, \dots, U_{n+1}, V_{n+1}$ be a feasible solution. 

As before, let $Z$ be our maximum rank slack in \eqref{p}. Since we can write $Z = B - \A x$ for some $x \in \rad{m},$ 
\eqref{eqn-ai0}  and \eqref{eqn-b0} imply
\begin{equation} \label{eqn-Z-UiVi=0} 
	\la Z, U_i + V_i \ra = 0 \; \text{for} \; i = 1,\ldots, n.
\end{equation}
Since $V_1 = 0, \, $ from \eqref{eqn-Z-UiVi=0} we deduce $\la Z, U_1 \ra = 0,$ hence the first $r$ rows and columns of $U_1$ are zero.
Let $Q_1$ be an invertible  matrix of suitably scaled eigenvectors of the lower right order $n-r$ block of $U_1, \, $
$$
T_1 = \bpx I_r & 0 \\
0 & Q_1 \epx,
$$
and apply the operations 
\begin{equation} \label{eqn-rescale-Ui} 
	\ba{rcl}
	U_i & := & T_1^{\top} U_i T_1 \\
	V_i  & := & T_1^{\top} V_i T_1 \; 
	\ena
\end{equation}
for all $i.$ Then we also apply the \eqref{eqn-rotate} operations 
with $T_1^{- \top}$ in place of $T.$ Afterwards  the $U_i$ and $V_i$ are still feasible in 
\eqref{d-ramana} and have the same objective value as they had before: this follows from 
the proof of Proposition \ref{prop-rotate-ramana}.

This rescaling does not change the maximum rank slack $Z$ in \eqref{p}. 
By Proposition \ref{prop-rotate-ramana}, it also does not change the optimal value of \eqref{d-ramana}. Finally, it does not change the value of 
\eqref{rel-d}: this follows by an argument similar to the one 
proving that rescaling  keeps the optimal value of \eqref{d} the same (this argument was given just 
before Proposition \ref{prop-rotate-ramana}). 

After this rescaling $U_1$ looks like 
\begin{equation}  \label{eqn-U1-rescaled} 
	U_1 = \bpx  0 & 0 \\
	0 & I_{r_1} 
	\epx, \; \text{where} \; 0 \leq  r_1 \leq n-r.
\end{equation}
Since $V_2 \in \tan(U_1), \,$ by part \eqref{tan-shape} in Proposition \ref{prop-basic} 
we deduce that only the last $r_1$ rows and columns of $V_2$ can be nonzero, hence $\la Z, V_2 \ra = 0.$ 
Again using \eqref{eqn-Z-UiVi=0} we deduce  $\la Z, U_2 \ra = 0.$ 

\myblue{
	Next we perform  the operations 
	\eqref{eqn-rescale-Ui} with  a suitable  invertible matrix $T_2$ in place of $T_1$ 
	to ensure $U_2 = 0 \oplus I_{r_2}$ for some $0 \leq r_2 \leq n-r$ 
	\footnote{Afterwards $U_1$ may not look like in equation \eqref{eqn-U1-rescaled} anymore, but for our purposes this does not matter.}. 
	We also apply the \eqref{eqn-rotate} operations 
	with $T_2^{- \top}$ in place of $T$ to keep the $U_i$ and $V_i$ feasible.  
}

Continuing, we produce matrices $T_3, \dots, T_n, \, $ and apply the operations \eqref{eqn-rescale-Ui} with $T_3, \dots, T_n$ in place 
of $T_1.$ We also apply the \eqref{eqn-rotate} operations with $T_3^{- \top}, \dots, T_n^{- \top} \, $
in place of $T.$ Afterwards, we have 
\begin{equation}
	U_n = \bpx  0 & 0 \\
	0 & I_{r_n} 
	\epx, \; \text{where} \; 0 \leq r_n \leq n-r.
\end{equation}
Since $V_{n+1}  \in \tan(U_n), $ again using part \eqref{tan-shape} in Proposition \ref{prop-basic} we deduce 
that only the last $r_n$ rows and columns of $V_{n+1}$ can be nonzero. Also, $U_{n+1} \succeq 0, $ so  
$U_{n+1} + V_{n+1}$ is a matrix whose upper left order $r$ block is psd, in other words 
$U_{n+1} + V_{n+1}$ is feasible in \eqref{rel-d}.

The proof is now complete.

\qed

Note that in Lemma \ref{lemma-valP<=valDram} we actually proved a stronger result, than what is strictly needed to prove 
\eqref{elaine}. Namely, we proved that after rescaling $\A$ and $B, \, $ and applying suitable similarity transformations to the $U_i$ and $V_i \, $
(to keep them feasible in \eqref{d-ramana}), the matrix $U_{n+1} + V_{n+1}$ is feasible in \eqref{rel-d}.

\myblue{
We can now prove the main result of the paper.
}

\pf{of Theorem \ref{theo:main}:} 
We have that 
\begin{equation}  \label{eqn-all-equal} 
	\val \eqref{p} \, = \, \val \eqref{rel-d} \, = \, \val \eqref{d-ramana},
\end{equation} 
where the first equality comes from Lemma \ref{lemma-relaxed-dual} and the second from Lemmas 
\ref{lemma-valP>=valDram} and \ref{lemma-valP<=valDram}.
We note that \eqref{eqn-all-equal} holds both when the optimal value of 
\eqref{p} is finite, and when it is $+ \infty.$ 

Further, if the optimal value of \eqref{p} is finite, then by Lemma 
\ref{lemma-relaxed-dual} 
the SDP  \eqref{rel-d} has a solution with the same value; and by Lemma 
\ref{lemma-valP>=valDram} the problem \eqref{d-ramana} has a solution with that value. This completes the proof.
\qed

In this work we assumed that \eqref{p} is feasible. 
On the other hand, when \eqref{p} is infeasible, Ramana's dual can provide a certificate to verify its infeasibility: for details we refer the reader to 
\cite{Ramana:97}.

 	\subsection{\bf SDP representation: proof of  Theorem \ref{theorem-repr}} 
 	\label{subsection-thm2} 
 
 Let us fix $U \succeq 0$ and let  $\T'(U)$ be the set on the right hand side of equation \eqref{eqn-UVW}. For 
 $W \in \rad{n \times n}$ and $\beta \in \rad{}, \,$ define the matrix 
 \begin{equation} \label{eqn-UWWtbeta} 
 M(W, \beta) \, := \, \bpx U & W \\ W^\top & \beta I \epx. 
 \end{equation}
  
 	To prove $\T(U) \supseteq \T'(U), \,$ suppose $W + W^\top \in \T'(U)$ and fix 
 	$\beta$ such that $M(W, \beta) \succeq 0.$
 
 		We want to show $\R(W) \subseteq \R(U), \, $ so to obtain a contradiction,  assume this is not the case. Then 
 the nullspace of $U$ is not contained in the nullspace of $W^\top, \, $ so we can choose $x$ such that 
 $Ux = 0$ and $W^\top x \neq 0.$ 
 
 Further, we pick $y$ such that $2 x^\top W y + \beta \norm{y}^2 < 0.$     
 Letting $z^\top := (x^\top, y^\top),$ we deduce
 $$
 z^\top M(W, \beta) z \, = \, 2 x^\top W y + \beta \norm{y}^2 < 0,
 $$
 the desired contradiction. 
   	
  	To show $\T(U) \subseteq \T'(U)$  suppose $W + W^\top \in \T(U), \, $ i.e., $\R(W) \subseteq \R(U).$ 
 	Hence $W = UH$ for some matrix $H.$ 
 	Define 
 	$$
 	R \, := \, \bpx I & - H \\
 	                      0 & I \epx.
 	                      $$
 	                      Then by an elementary calculation,   
 	 	                      \begin{equation} \nonumber 
 	                      R^\top M(W, \beta) R \, = \, \bpx U & 0 \\ 0 & \beta I - H^\top U H \epx,
 	                      \end{equation}
 	so if $\beta$ is large enough, then $\beta I - H^\top U H \succeq 0. \, $ 
 	Hence 
 	$M(W, \beta) \succeq 0, \, $ so $W + W^\top \in \T'(U), \, $ 
 	completing the proof.
 	\qed

 	\section{Conclusion: a larger example, and the tangent space} 
 	\label{section-conclude} 
 	
 	We first present a larger example to illustrate Ramana's dual and walk the reader through all steps in our previous discussions. 
 	
 		\bex \label{ex-larger}  Consider the SDP 
 		\beq    \label{problem-gap-n=3} 
 		\ba{rl} 
 		\sup &  x_2 +  x_3 \\
 		s.t. & x_1 \underbrace{\bpx 1 & \pha{0}  & \pha{\ti} & \pha{0} \\
 		\pha{0} & 0  & \pha{0} & \pha{0} \\
 		\pha{0} & \pha{0} & 0  & \pha{0} \\
 		\pha{0} & \pha{0} & \pha{0} & 0  \epx}_{A_1}
 		+ x_2 \underbrace{\bpx 0  & \pha{0} & 1   & \pha{0}   \\
 		\pha{0} & 1 & \pha{0} & \pha{0} \\
 		1  & \pha{0} & {0}  & \pha{0} \\
 		\pha{0}   & \pha{0}  & \pha{0} & {0}     \epx}_{A_2} 
 		+ x_3 \underbrace{\bpx 0  & \pha{0} & \pha{0} & \pha{0} \\
 		\pha{0} & 0 & \pha{0} & 1 \\
 		\pha{0} & \pha{0} & 1 & \pha{0} \\
 		\pha{0} & 1 & \pha{0} & 0 \epx}_{A_3} \preceq 
 		\underbrace{\bpx 1 & \pha{0} & \pha{0} & \pha{0} \\
 		\pha{0} & 1 & \pha{0} & \pha{0} \\
 		\pha{0} & \pha{0} & 0 & \pha{0} \\
 		\pha{0} & \pha{0} & \pha{0} & 0  \epx}_{B}.  
 		\ena 
 		\eeq 
 		We will proceed as follows. We first show that 	this SDP has a positive duality gap. Then we calculate the key players in the paper: the maximum rank slack of 
 		\eqref{problem-gap-n=3} in the form given in Assumption \ref{assumption-slack}; an optimal solution of the strong dual \eqref{rel-d};   the $Y_i$ certificates for the maximum rank slack (given in Lemma \ref{lemma-complete-certificate}); and finally, an optimal solution to Ramana's dual.
 		\begin{enumerate}
 			\item \label{lastex-item1} 
 			To calculate the primal optimal value we let 
 		$x$ be feasible in \eqref{problem-gap-n=3}, and 
 			\begin{equation} \label{eqn-S} 
 				S \, := \, \bpx 1 - x_1 & \pha{0} & - x_2  & \pha{0}  \\
 			\pha{0} & 1- x_2  & \pha{0} & - x_3  \\
 			- x_2   & \pha{0} &  - x_3  & \pha{0} \\
 			\pha{0}   & - x_3  & \pha{0} & 0  \epx 
 				\end{equation}
 			the corresponding slack matrix. 	Since the lower right corner of $S$ is $0, \, $ and $S \succeq 0, \, $ 
 			the last row and column of $S$ is zero. 
 			Thus $x_3=0, \, $ so the $(3,3)$ element of $S$ is $0. $ Hence the third row and column of $S$ are also zero, so $x_2 = 0.$ 
 			
 			We have learned that in any feasible solution 
 			\begin{equation} \label{eqn-x2x3=0} 
 				x_2 = x_3 = 0, \,
 			\end{equation}
 			 so the objective function is identically zero on the primal feasible set.
 			\item 
 			To calculate the dual optimal value,  suppose $Y \in \psd{4}$ is feasible in the dual. 
 			Since  $\la A_1, Y \ra = y_{11}=0, $ and  $Y \succeq 0, \, $
 			the first row (and column) of 
 			$Y$ is zero. Thus, we have $\la A_2, Y \ra = y_{22} = 1, $ so the dual objective function is identically $1$ on the feasible set. 
 			
 			For example, 
 			\begin{equation} \label{eqn-Y-large} 
 				Y \, = \, \bpx {0}  & \pha{0}  & \pha{0}   & \pha{0}   \\
 				\pha{0}   & 1  & \pha{0} & \pha{0}    \\
 				\pha{0}  & \pha{0} &  1    & \pha{0} \\
 				\pha{0}   &  \pha{0}  & \pha{0} & 0  \epx
 			\end{equation}
 			is an optimal solution in the dual. 

 		\item \myblue{Next we compute the optimal value of \eqref{rel-d}. Because of \eqref{eqn-x2x3=0},  the maximum rank slack in the SDP 
 			\eqref{problem-gap-n=3} is just the right hand side.}
 	 	Thus, in the strong dual \eqref{rel-d} only the upper left $2 \times 2$ block of $Y$ must be positive semidefinite. It follows that 
 		\begin{equation}
 			Y^*  \! := \!  \bpx 0 & 0 & 0.5  & 0  \\
 			0  & 0  & \pha{0} & \pha{0}    \\
 			0.5   & \pha{0} &  1    & \pha{0} \\
 			{0}  &  \pha{0}  & \pha{0} & 0  \epx
 		\end{equation}
 		is optimal in \eqref{rel-d} with value zero. 
 		
 	\item The $Y_i$ matrices from Lemma \ref{lemma-complete-certificate} below certify that the right hand side in \eqref{problem-gap-n=3} is the maximum rank slack:
 	\begin{equation}
 		\ba{rclrcl}
 		Y_1  & \! := \! & \bpx 0 & 0 & 0  & 0  \\
 		0  & 0  & \pha{0} & \pha{0}    \\
 		0   & \pha{0} &  0    & \pha{0} \\
 		{0}  &  \pha{0}  & \pha{0} & 1  \epx,  \\
 		Y_2  & \! := \! & \bpx 0 & 0 & 0  & 0  \\
 		0  & 0  & \pha{0} & -1     \\
 		0   & \pha{0} &  2    & \pha{0} \\
 		{0}  &  -1   & \pha{0} & 0  \epx.  \\
 		\ena 
 	\end{equation}
 Indeed, the $Y_i$ are of the form required in \eqref{eqn-Yi} and satisfy the equations  \eqref{eqn-Yi-A*} 
 \footnote{The following argument may better explain the role of the $Y_i$ matrices. 
 	If $S$ is any slack, then 	$\la S, Y_1 \ra = \la S, Y_2 \ra = 0.  \, $ We invite the reader to check that  these 
 	equations lead to the same argument that we gave 
 	in paragraph 	\eqref{lastex-item1} that show the last two rows and columns of $S$ are zero.}.

\item We finally construct a solution of \eqref{d-ramana} with value $0,$ following the proof of Lemma \ref{lemma-valP>=valDram}. 
We first set 
	$U_i$ and $V_i$ to $0$ for $i=0,1,2. $ 
	
	Then we split the $Y_i$ and $Y^*$ to define the other $U_i$ and $V_i$ in \eqref{d-ramana}: 
		\begin{equation}
		\ba{rclcl}
		Y_1 & = &  \underbrace{\bpx 0 & 0 & 0  & 0  \\
		0  & 0  & \pha{0} & \pha{0}    \\
		0   & \pha{0} &  0    & \pha{0} \\
		{0}  &  \pha{0}  & \pha{0} & 1  \epx}_{U_3}  & + & \underbrace{0}_{V_3},  \\ 
		Y_2   & \! := \! & \underbrace{\bpx 0 & 0 & 0  & 0  \\
		0  & 0  & \pha{0} & \pha{0}    \\
		0   & \pha{0} &  2    & \pha{0} \\
		{0}  &  \pha{0}  & \pha{0} & 0  \epx}_{U_4}  & + &  \underbrace{\bpx 0 & 0 & 0  & 0  \\
		0  & 0  & \pha{0} & -1     \\
		0   & \pha{0} &  0    & \pha{0} \\
		{0}  & -1   & \pha{0} & 0   \epx}_{V_4},  \\ 
		Y^* & \! := \! & \underbrace{\bpx 0 & 0 & 0  & 0  \\
		0  & 0  & \pha{0} & \pha{0}    \\
		0   & \pha{0} &  1    & \pha{0} \\
		{0}  &  \pha{0}  & \pha{0} & 0  \epx}_{U_5}   & + & \underbrace{ \bpx 0 & 0 & 0.5  & 0  \\
		0  & 0  & \pha{0} & \pha{0}      \\
		0.5   & \pha{0} &  0    & \pha{0} \\
		{0}  & \pha{0}    & \pha{0} & 0   \epx}_{V_5}.
		\ena 
	\end{equation}
Note that by part \eqref{tan-shape} in Proposition \ref{prop-basic} we have 
$V_4 \in \tan(U_3)$ and $V_5 \in \tan(U_4), $ thus 
$U_0, V_0, \dots, U_5, V_5$ is indeed a solution of \eqref{d-ramana} with value $0. \, $ 

\end{enumerate}
 		
 	\end{Example}

 		\myblue{
 			In the last remark we describe the geometry of the  mysterious looking set   $\tan(U)$
 			that plays a crucial role in Ramana's dual. }
 		
 				Suppose  $U \in \psdn.$ 
 				Recall that we call $\tan(U)$ the {\em tangent space of $\psd{n}$ at $U$}  and we gave a 
 			purely algebraic definition  in  \eqref{eqn-tan:2}. An equivalent geometric expression is 
 				\begin{equation} \label{eqn-tandef-geo} 
 			\T(U) \, = \,  \Bigl\{ \, V \in \symn \, \Bigl| \, \dfrac{1}{\epsilon} \dist( U \pm \eps V, \psdn)   \rightarrow 0 \, \text{as} \, \eps \searrow 0   \, \Bigr\},
 				\end{equation}	
 				where $\dist(X, \psd{n}) \, = \, \inf  \{ \norm{X - Y} \, | \, Y \in \psd{n} \, \}$ is the distance 
 				of matrix $X$ from $\psd{n}$ \footnote{We can use any matrix norm, for example the spectral norm or the Frobenius norm.}. For a detailed proof of the equivalence, see e.g. \cite[Lemma 3]{Pataki:17}. 
 				
 				For example, if 
 				$$
 				U = \bpx 0 & 0 \\ 0 & 1 \epx, \, V = \bpx 0 & 1 \\ 1 & 0 \epx, \, 
 				$$
 				then $V \in \T(U) \, $ according to the algebraic definition  \eqref{eqn-tan:2}. We also have 
 				 $V \in \T(U) \, $ according to the geometric definition \eqref{eqn-tandef-geo}, since we only need to change the upper left corner of $U \pm \epsilon V$ to $\epsilon^2$ to make it psd.

\begin{figure}[H]
 						\centering
 						\includegraphics[width = 6cm, height=5cm]{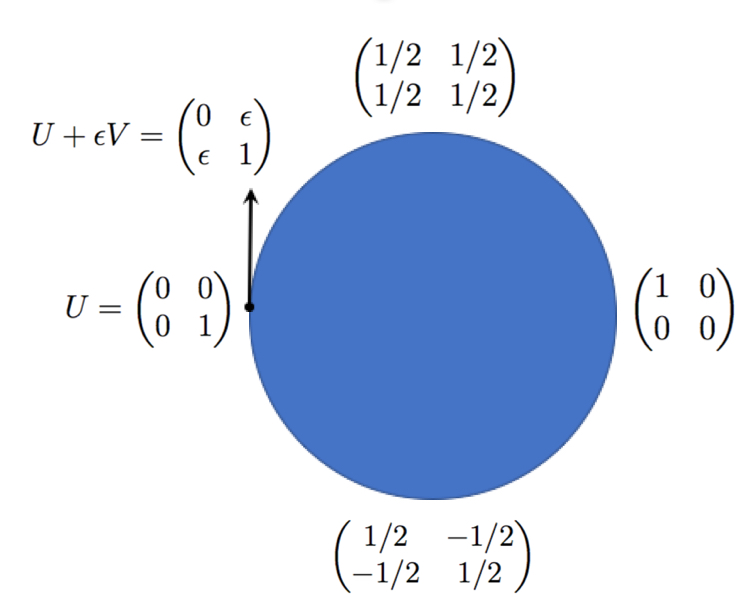} 
 						\caption{The matrix $U + \epsilon V$ is ``almost" psd, but not quite}
 						\label{figure-tangent}  
 					\end{figure}

 				We illustrate this on Figure \ref{figure-tangent}. We let $C$ be the set of psd matrices with trace $1, \, $ and we describe  
 				$C$ with just two parameters,  as 
 				$$
 				C \, = \, \Biggl\{ \, \bpx x & y \\ y & 1 - x \epx  \, \mid \, 1 \geq x \geq 0, \, x(1-x) - y^2 \geq 0 \, \Biggr\}.  
 				$$
 				Since we can rewrite the quadratic inequality as $( x - \frac{1}{2})^2 + y^2 \leq \frac{1}{4}, \, $ the set $C$ 
 				is a circle of radius $\frac{1}{2}$ centered at $(\frac{1}{2}, 0).$ Figure \ref{figure-tangent} shows $C$ together with $U$ and $U + \epsilon V$ for a small $\epsilon > 0.$ 

 				Nevertheless, it is interesting that one can completely derive Ramana's dual using only purely linear algebraic arguments about $\T(U), \, $ rather than any geometric ones, and this is the route we followed.

 			\section*{Acknowledgements}
 			Part of this paper was written during a visit to Chapel Hill by the first author and he would like to express his gratitude to
 			Prof.~Takashi~Tsuchiya, for creating conditions that made this visit possible.
 			The  work of the first author was partially supported by the JSPS KAKENHI Grant Numbers JP15H02968 and JP19K20217.
 			
 			The second author thanks Pravesh Kothari and Ryan O' Donnell for helpful discussions on SDP.
 			The work of the second author was supported by the National Science Foundation, award DMS-1817272.
 			
 			Both authors are grateful to Siyuan Chen, Alex Touzov,  and Yuzixuan Zhu for their careful reading of the manuscript and their helpful comments.
 			
 			Most importantly, we are very grateful to the anonymous referees whose comments  and suggestions greatly improved the manuscript. 
 		\bibliographystyle{plain}
 		\bibliography{mysdpMelody}
 		 
 \end{document}